\documentclass[12pt]{article}

\usepackage{amsmath, amssymb, amsthm}

\theoremstyle{definition}

\numberwithin{equation}{section}

\newtheorem{thm}{THEOREM}[section]

\newtheorem{lem}[thm]{Lemma}

\newtheorem{prop}[thm]{PROPOSITION}

\theoremstyle{definition}
\newtheorem{defn}{Definition}[section]
\theoremstyle{remark}
\newtheorem{rem}{Remark}[section]

\newcommand{\tref}[1]{Theorem~\ref{#1}}
\newcommand{\cref}[1]{Korollar~\ref{#1}}

\newcommand{\lref}[1]{Lemma~\ref{#1}}

\title{$C^*$-algebras coming from some buildings.}
\author{{\sc  Alina Vdovina} \\
\small  {Mathematisches Institut} \\
\small  {Beringstrasse 1, 53115 Bonn} \\
\small  { e-mail: \ \ alina@math.uni-bonn.de}
}
\date{}

\begin{document}

\maketitle

\abstract{
We construct compact polyhedra with
triangular faces whose 
links are generalized $3$-gons.
They are interesting
compact spaces covered by Euclidean buildings of type $A_2$.
Those spaces give us two-dimensional subshifts,
which can be used to construct some $C^*$-algebras.

\section{Introduction}

 Given a graph $G$ we assign to each edge the length $1$. The
diameter of the graph is its diameter as a length metric space,
its  injectivity radius is half of the length of the smallest circuit.

The following 
definition is equivalent to the usual one.

\begin{defn}
 For a natural number $m$ we call a connected graph $G$ a generalized
$m$-gon, if its diameter and injectivity radius are both equal to $m$.
\end{defn}

 A graph is {\em bipartite} if its set of vertices
can be partitioned into two disjoint subsets $P$ and $L$ such that no
two vertices in the same subset lie on a common edge.
 Such a graph can be interpreted as a planar geometry, i.e.
a set of points $P$ and a set of lines $L$ and an incidence relation
$R\subset P\times L$. On the other hand each planar geometry
can be considered as a bipartite graph.

 Under this correspondence projective planes are the same as
generalized $3$-gones.

 Let $G$ be a planar geometry.
 For a line $y\in L$ we denote by $I(y)$ the set of all
points $x\in P$ incident to $y$. If no confusion can arise
we shall write $x\in y$ instead of $x\in I(y)$ and $y_1\cap
y_2$ instead of $I(y_1)\cap I(y_2)$.
 A subset $S$ of $P$ is called collinear if it is contained
in some set $I(y)$, i.e. if all points of $S$ are incident to a line.

 Given a planar geometry $G$ we shall denote by $G^{\prime}$ its
dual geometry arising by calling lines resp. points of $G$ points
resp. lines of $G^{\prime}$. The  graphs corresponding to
$G$ and  $G^{\prime}$ are isomorphic.

We will call a {\em polyhedron} a two-dimensional
 complex which is obtained from several oriented $p$-gons
by identification of corresponding sides.
Consider a point of the polyhedron and 
take a sphere of a small radius at this point.
The intersection of the sphere with the polyhedron is
a graph, which is called the {\em link} at this point.

\bigskip

\noindent
{\bf Definition.} Let $\tilde A_2$ be a tessellation of the
Euclidean plane by 
regular triangles.
A {\em Euclidean  building} of type  $\tilde A_2$
is a polygonal complex $X$,
which can be expressed as the union of subcomplexes called 
apartments such that:

1. Every apartment is isomorphic to $\tilde A_2$.

2. For any two polygons of $X$, there is an apartment
containing both of them.

3. For any two apartments $A_1, A_2 \in X$ containing
the same polygon, there exists an isomorphism $ A_1 \to A_2$
fixing $A_1 \cap A_2$.

\medbreak

If we consider a polyhedron with triangular
faces and incidence graphs of finite projective
planes as links, the universal covering of the polyhedron
 is an Euclidean building $\Delta$, see \cite{BB}, \cite{Ba}.

So, to construct Euclidean buildings
with compact quotients, it is sufficient to construct
finite polyhedra with appropriate links.

We construct a family of 
compact polyhedra
with $3$-gonal faces whose links
are generalized $3$-gons.
Then we will use those polyhedra to construct
systems of 2-dimensional words, which give examples of new $C^*$-
algebras.

\bigskip

We recall the definition of the polygonal presentation, given in \cite{V}.

\medskip

\noindent
{\bf Definition.} Suppose we have $n$ disjoint connected bipartite graphs
 $G_1, G_2, \ldots G_n$.
Let $P_i$ and $L_i$ be the sets of black and white vertices respectively in
$G_i$, $i=1,...,n$; let $P=\cup P_i, L=\cup L_i$, $P_i \cap P_j = \emptyset$
 $L_i \cap L_j = \emptyset$
for $i \neq j$ and
let $\lambda$ be a  bijection $\lambda: P\to L$.

A set $\mathcal{K}$ of $k$-tuples $(x_1,x_2, \ldots, x_k)$, $x_i \in P$,
will be called a {\em polygonal presentation} over $P$ compatible
with $\lambda$ if

\begin{itemize}

\item[(1)] $(x_1,x_2,x_3, \ldots ,x_k) \in \mathcal{K}$ implies that
   $(x_2,x_3,\ldots,x_k,x_1) \in \mathcal{K}$;

\item[(2)] given $x_1,x_2 \in P$, then $(x_1,x_2,x_3, \ldots,x_k) \in 
\mathcal{K}$
for some $x_3,\ldots,x_k$ if and only if $x_2$ and $\lambda(x_1)$
are incident in some $G_i$;

\item[(3)] given $x_1,x_2 \in P$, then  $(x_1,x_2,x_3, \ldots ,x_k) 
\in \mathcal{K}$
    for at most one $x_3 \in P$.

\end{itemize}

If there exists such $\mathcal{K}$, we will call $\lambda$ a 
{\em basic bijection}.

Polygonal presentations for $n=1$, $k=3$ were listed in $\cite{Cart}$
with the incidence graph of the finite projective plane of order two 
or three as the graph $G_1$.

One can associate  a polyhedron $X$ on $n$ vertices with
each polygonal presentation $\mathcal{K}$ as follows:
for every cyclic $k$-tuple $(x_1,x_2,x_3,\ldots,x_k)$ from
the definition
we take an oriented $k$-gon on the boundary of which
the word $x_1 x_2 x_3\ldots x_k$ is written. To obtain
the polyhedron we identify the sides with the same label of our
polygons, respecting orientation.
We will say that the
polyhedron $X$ {\em corresponds} to the polygonal
presentation $\mathcal{K}$.

The following lemma was proved in \cite{V}:
\begin{lem} \label{Main}
 A polyhedron $X$ which corresponds to
a polygonal presentation $\mathcal{K}$ has
  graphs $G_1, G_2, \ldots, G_n$ as the links.
\end{lem}

\noindent
{\bf Remark.} Consider a  polygonal
presentation $\mathcal{K}$. Let $s_i$ be the number of vertices
of the graph $G_i$ and $t_i$ be the number of edges of $G_i$,
$i=1,...,n$.
If the polyhedron $X$  corresponds to the polygonal
presentation $\mathcal{K}$, then $X$ has $n$ vertices
(the number of vertices of $X$ is equal to the number of graphs),
$k \sum_{i=1}^n s_i$ edges and $\sum_{i=1}^n t_i$ 
faces, all faces are polygons with $k$ sides.

 \section{Construction of the polyhedra.}

 Let $G$ be a finite projective plane and let $P$ resp. $L$ denote the
set of  its points resp. lines.

   Assume that a bijection $T:P\to L$ is given and satisfies the following
properties 
\begin{enumerate}
\item For each $x\in P$ the point $x$ and the line $T(x)$ are not incident.

\item For each pair $x_1,x_2$ of different points  in $P$ the points
$x_1,x_2 $ and $T(x_1)\cap T(x_2)$ are not collinear.   
\end{enumerate}

\medskip

Projective planes with such bijections exist, see for example \cite{V2}.

\medskip

\begin{lem} \label{easy}
 Let $T:P\to L$ be as above, $y\in L$ a line. Then the map 
$T^*:I(y) \to I(y)$ given by $T^*(x)=T(x)\cap I(y)$ is a
bijection.  
\end{lem}
\begin{proof}
 By the first property of $T$ the map $T^*$ is well defined, 
by the second property it  must be injective. Since $I(y)$ is
finite, the statement follows.  
\end{proof}

 Let $G,P,L, T:P\to L$ be as above. Let   $P =\{x_1,..x_p \} $ be a labeling 
of points in $P$  and set $y_i=T(x_i)$.
 Consider the following
set  $O\subset P\times P \times P$, consisting of all
triples $(x_i,x_j,x_k)$ satisfying  $x_i\in y_k$, $x_j\in y_i$ and
 and $x_j \in y_k$.

\begin{rem}
 The conditions on $(x_i,x_j,x_k)\in K$ are not cyclic.
We require $x_j \in y_k$ and not $x_k \in y_j$ !!  For this
reason in the polygonal presentations defined below dual 
graphs of $G$ appear. 
\end{rem}

The following lemma is crucial for the later construction:
\begin{lem} \label{crucial}
 A pair $(x_i,x_k)$ resp. $(x_i,x_j)$ resp. $(x_j,x_k)$ is a part
of at most one  triple $(x_i,x_j,x_k)\in K$ and such a triple exists
iff $x_i \in y_k$ resp. $x_j \in y_i$ resp. $x_j \in y_k$ holds.  
\end{lem}

\begin{proof}
 The conditions stated at the end are certainly necessary. 

1) Let $x_i\in y_k$ be given. Then $y_i$ and $y_k$ are different and the
point $x_j  =y_i\cap y_k$ is uniquely defined. 

2) Let $x_j \in y_i$ be given. Then $x_j$ and $x_i$ are different, so
there is exactly one line $y_k$ containing $x_j$ and $x_l$.

3) Let $x_j \in y_k$ be given. Then $(x_i,x_j,x_k)$ is in $K$ iff
for the map $T^*:I(y_k)\to I(y_k)$  of \lref{easy} 
the equality $T^* (x_i)=x_j $ holds. By \lref{easy}  the point
$x_i$ is uniquely defined. 
\end{proof}

 Now we are ready for the polygonal presentations.
Let the notations be as above, $G_1$ and $G_2$ two projective planes with
isomorphism $J^t:G\to G_t$ and $G_3$ a projective plane with an
isomorphism $J^3:G' \to G_3$ of the dual projective plane $G'$ of $G$.
For $t=1,2$ we set $x^t_i =J^t(x_i)$, $y^t_i=J^t(y_i)$ and for $t=3$
we set $x^3_i =J^3 (y_i)$ and $y^3_i =J^3(x_i)$.

 Let $P_t$ resp. $L_t$ be the set of lines of $G_t$. For
$P=\cup P_t$ and $L=\cup L_t$ we consider the bijection
$\lambda :P\to L$ given by  $\lambda (x_i^t) = y_i^{t+1}$
 ($t+1$ is taken modulo $3$).

 Now consider the subset  $\mathcal{T}$ of $P\times P \times P$ consisting
of all triples $(x^1_i,x^2_j,x^3_k)$ with $(x_i,x_j,x_k)\in K$ and
all cyclic permutation of such triples.

 The statement of \lref{crucial} can be now reformulated as:

\begin{prop}
 The subset  $\mathcal{T}$ of $P\times P \times P$ defines a
polygonal presentation compatible with $\lambda$.
\end{prop}

The polyhedron $X$ which corresponds to 
$\mathcal{T}$ by the construction of \lref{Main} has 
triangular faces and exactly  three vertices with two links
naturally isomorphic to $G$ and one link naturally isomorphic
to the dual $G'$ of $G$.  By \cite{BB},\cite{Ba} the universal 
covering of $X$  is a 
Euclidean  building of type $\Tilde{A}_2$.

\section{Subshift coming from a polygonal presentation.}
Let  $\mathcal{T}$ be a polygonal presentation with $n=3$, $k=3$,
where all there graphs $G_1$, $G_2$ and $G_3$ are incidence
graphs of finite projective planes of order $q$. 
The polyhedron, which corresponds
to $\mathcal{T}$, has triangular faces and three vertices.
We will consider polyhedra such that all three vertices of each 
triangle have different graphs as links.
In this case  
we can give a Euclidean metric to every face. In this metric
 all sides of the triangles
are geodesics of the same length. The universal covering of the polyhedron
 is an Euclidean building $\Delta$, see \cite{BB}, \cite{Ba}.
 Each element of $\mathcal{T}$ may be identified with  
an oriented basepointed triangle in $\Delta$. We now construct
a 2-dimensional shift system associated with  $\mathcal{T}$.
The transition matrices $M_1, M_2$
in the way, defined as in \cite{RS2}: if $\alpha=(x_1,x_2,x_3),
 \beta=(y_1, y_2, y_3) \in \mathcal{T}$
say that $M_1(\alpha, \beta)=1$ if and only if there exists
  $\psi=(x_2,y_1,z_3)$ and $M_1(\alpha, \beta)=0$ otherwise.
 In a similar way,
$M_2(\alpha,\gamma)=1$ for $\alpha=(x_1,x_2,x_3), \gamma=(z_1,z_2,z_3)$
if and only if there exists $\psi=(x_2,y_1,z_3)$ and
$M_2(\alpha,\gamma)=0$
otherwise. The matrices  $M_1, M_2$ of order
$3(q+1)(q^2+q+1) \times 3(q+1)(q^2+q+1)$ are nonzero
$\{0,1\}$ matrices.
We will use $\mathcal{T}$ as an alphabet and $M_1, M_2$
as transition matrices to build up 2-dimensional words as in
\cite{RS1}. Let $[m,n]$ denote $\{m, m+1,...,n\}$, where $m \leq n$
are integers. If $m,n \in \mathbb{Z}^2$, say that $m \leq n$
if  $m_j \leq n_j$ for j=1,2, and when $m \leq n$ let
$[m,n]=[m_1,n_1] \times [m_2,n_2]$. In $ \mathbb{Z}^2$, let 0
denote the zero vector and let $e_j$ denote the $j$-th standard 
unit basis vector. If $m \in \mathbb{Z}_+^2=\{m \in \mathbb{Z}^2;
m \geq 0\}$, let

\noindent $W_m=\{w:[0,m] \to \mathcal{T}; M_j(w(l+e_j),w(l)=1$ where
$l,l+e_j \in [0,m]\}$
and call the elements of $W_m$ words.

 In order to apply the theory
from \cite{RS1} we need the matrices $M_1,M_2$ to satisfy the following
conditions:

\vfill\eject

(H0) Each $M_i$ is a nonzero $\{0,1\}$-matrix.

(H1a) $M_1M_2=M_2M_1$.

(H1b) $M_1M_2$ is a $\{0,1\}$-matrix.

(H2) The directed graph with vertices $\alpha \in \mathcal{T}$
and directed edges $(\alpha,\beta)$
whenever $M_i(\alpha, \beta)=1$ for some $i$ is irreducible.

(H3) For any nonzero $p \in \mathbb{Z}^2$, there exists
a word $w \in W$ which is not $p-periodic$, i.e., there exists $l$
so that $w(l)$ and $w(l+p)$ are both defined but not equal.

\bigbreak

In \cite{RS1} some $C^*$-algebra is defined by
 the system of words $W_m$, where $m \in \mathbb{Z}_+^2$.
It is proved there, that if the matrices  
$M_1,M_2$  satisfy the conditions
(H0),(H1a,b),(H2),(H3), then this algebra is simple,
purely infinite and nuclear.

\medbreak

Now we prove the conditions (H0),(H1a,b),(H2),(H3)
for our two-dimensional shift.
By definition of matrices $M_1,M_2$ they are nonzero
$\{0,1\}$ matrices, so (H0) holds. 
If we have $\alpha$, $\beta$, $\psi$, such that  $M_1(\alpha, \beta)=1$,
$M_2(\beta, \psi)=1$, then $\gamma$ such 
that $M_2(\alpha, \gamma)=1$, $M_1(\gamma, \psi)=1$, is uniquely
defined because of properties
of finite projective planes. Conditions (H1a,b) follow.
To prove (H2) we need to color sides of triangles in three
different colors. This is possible since there are three vertices
in the polyhedron with different graphs as links.
So, all triangles from $\mathcal{T}$ have one of three possible
colorings.
We need to show, that for any $\alpha, \beta \in \mathcal{T}$
we can choose $r>0$ such that $M_j^r(\alpha, \beta)>0$, where $j=1,2$.
Geometrically it means that any  $\alpha, \beta \in \mathcal{T}$
can be realized so that $\beta$ lies in some sector with base $\alpha$
(for more details see \cite{RS1}). Without loss of generality
we can assume, that $j=1$. We will say, that $\beta \in \mathcal{T}$
is reachable from $\alpha \in \mathcal{T}$ in $r$ steps, if there is
 $r>0$ such that $M_1^r(\alpha, \beta)>0$. It is easy to see, that
every triangle is reachable from some triangle of other color
in one or two steps. So, to prove (H2) we need to show, that
any triangle is reachable from another one of the same color.
Now we can use the proof of the Theorem 1.3 from \cite{RS1},
since at each step of this proof it is only used, that
the link at each vertex of the building is an incidence graph
of a finite projective plane, which is true in our case too.
The proof of (H3) is identical to the proof of (H3)
in the case of the subshift considered in \cite{RS2}.

\end{document}